\documentclass[final,onefignum,onetabnum]{siamart190516}

\usepackage[a4paper,top=2.0cm,bottom=2.0cm,left=2.3cm,right=2.0cm]{geometry}
\usepackage{lipsum}
\usepackage{amsfonts,physics}
\usepackage{graphicx,subfigure}
\usepackage{epstopdf}
\usepackage{algorithmic}  
\usepackage{tikz}
\usepackage{mathtools}
\usepackage{booktabs}
\usepackage{multirow}
\usepackage{siunitx}
\usetikzlibrary{calc}
\usetikzlibrary{backgrounds}
\usepackage{pgfplots,pgfplotstable}
\usetikzlibrary{spy}  
\pgfplotsset{compat=1.5.1}
\pgfplotsset{every axis/.append style={
                    xlabel={$x$},          
                    ylabel={$y$},          
                    label style={font=\small},
                    title style={font=\small},
                    tick label style={font=\tiny}  
                    }}

\usepackage{newfloat}
\ifpdf
  \DeclareGraphicsExtensions{.eps,.pdf,.png,.jpg}
\else
  \DeclareGraphicsExtensions{.eps}
\fi

\newsiamremark{remark}{Remark}
\newsiamremark{hypothesis}{Hypothesis}
\crefname{hypothesis}{Hypothesis}{Hypotheses}
\newsiamthm{claim}{Claim}

\usepackage[framemethod=TikZ]{mdframed}
\newcounter{exam}[section] 
\renewcommand{\theexam}{\arabic{exam}}
\numberwithin{exam}{section}
\crefname{exam}{example}{examples}
\Crefname{exam}{Example}{Examples}

\creflabelformat{exam}{(#2#1#3)}
\newenvironment{exam}[2][]{%
\refstepcounter{exam}%
\ifstrempty{#1}%
{\mdfsetup{%
frametitle={%
\tikz[baseline=(current bounding box.east),outer sep=0pt]
\node[line width=2pt,anchor=east,rectangle,draw=blue!20,fill=blue!20]
{\strut Example~\theexam};}}
}%
{\mdfsetup{%
frametitle={%
\tikz[baseline=(current bounding box.east),outer sep=0pt]
\node[line width=2pt,anchor=east,rectangle,draw=blue!20,fill=blue!20]
{\strut Example~\theexam:~#1};}}%
}%
\mdfsetup{innertopmargin=10pt,linecolor=blue!20,%
linewidth=2pt,topline=true,%
frametitleaboveskip=\dimexpr-\ht\strutbox\relax
}
\begin{mdframed}[]\relax%
\label{#2}}{\end{mdframed}}

\headers{Time-continuous dual estimation}{Angwenyi David}

\title{Feedback Particle Filter With Stochastically Perturbed Innovation And Its Application to Dual Estimation \thanks{Submitted to the editors on \today.
\funding{This work was funded by the DAAD and the Kenya National Research Fund (NRF)}}}

\author{Angwenyi David\thanks{Masinde Muliro University of Science and Technology 
  (\email{dangwenyi@mmust.ac.ke}).}}

\usepackage{amsopn}

\ifpdf
\hypersetup{
  pdftitle={Comparison of the performance of different filters in time-continuous dual estimation},
  pdfauthor={Angwenyi David}
}
\fi

\begin{document}

\maketitle

\begin{abstract}
Particle filters have, in recent years, been found to perform well in highly nonlinear problems as well as in estimation of parameters. However, there is still the problem of particle degeneracy in particle filters which has led to the invention of, among others, feedback particle filters. In this paper, we introduce a stochastically perturbed feedback particle filter and show that it is exact. The novelty is in the fact that the innovation process is stochastically perturbed. Resampled sinkhorn particle filter is also introduced. We then compare their performance with that of other filters in simultaneous state and parameter estimation. 
\end{abstract}



\section{Introduction}
Most state-space models are characterized by, among other things, parameters---which can be constant or varying. By a parameter is comprehended and signified a measurable factor which defines a model and influences its operation. As the parameter changes, so does the model; expressed differently, a parameter is unique to a model it characterises. The choice of a certain model, therefore, is achieved by choosing the right parameters. It occurs more often than not---in hidden Markov models, for instance---that measurements are available but the underlying signal is not readily apparent. This forms an example where parameter estimation is paramount: measurements are used to learn the model parameters, which, in turn, are used to fit the model. 

For the ensuing discussion, we modify the usual signal and measurement equations, to include parameters, now and henceforth signified by a $d-$dimensional vector, $\theta$, thus
\begin{subequations}
\begin{align}
\textrm{\textcolor{black}{Signal}:} & \quad dx_{t} = f(x_{t}, \theta)dt  + g(x_{t})d\beta_{t}; \quad t_{0} \leq t, \label{eq:5.1.1a}\\
\textrm{\textcolor{black}{Measurement}:}& \quad dy_{t} = h(x_{t}, \theta)dt + R^{1/2}(t)d\eta_{t};  \quad t_{0} \leq t, \label{eq:5.1.1b}
\end{align}
\end{subequations}
where
\begin{center}
\begin{tabular}{lll}
\textbf{Term} & \textbf{Name} & \textbf{Dimension}\\ 
$x_{t}$ & state vector & $n \times 1$\\ 
$f(x_{t},t)$ & drift function & $n \times 1$ \\ 
$g(x_{t},t)$ & diffusion function & $n \times m$ \\ 
$\{\beta_{t},\, t>t_{0}\}$ & Brownian motion process & $m \times 1$ \\ 
$y_{t}$ & output vector & $r \times 1$ \\ 
$h(x_{t},t)$ & sensor function & $ r \times 1$ \\ 
$R(t)$ & time-function matrix & $r \times r$ \\
$\{\eta_{t},\, t>t_{0}\}$ & standard Brownian motion process & $r \times 1$\\ 
\end{tabular} 
\end{center}

Parameter estimation problem concerns finding the optimal parameter so that the signal best fits the data \cite{Ein,Sar}. This is, classically, achieved by optimization procedure where a cost function is minimized \cite{Lew}. The cost function mostly defines the discrepancy between the state and the measurements. Intuitively, parameter estimation can be seen as a procedure for seeking a parameter value that gives the least discrepancy between the state and the corresponding measurements (also known as the algebraic distance or the residual). The method of least squares has been extensively used to define an objective function. Given the increment in measurement, $dy_{t}$, of the state, $x_{t}$, at time $t$, the objective function, $\mathcal{J}(\theta)$, in the least squares sense, is given by
\begin{equation}
 \mathcal{J}(\theta) = \int_{t_{0}}^{t}w_{t}\Vert dy_{t} - h(x_{t}, \theta ) dt\Vert^{2},
 \end{equation} 
where $w_{t}$ is the weighting function.

Most commonly used procedures in the framework of least-squares include: linear least-squares, orthogonal least-squares, gradient weighted least-squares, bias corrected renormalization. This paper, however, attends not to the study of least-squares approaches. Suffice it to only direct the interested reader to the article \cite{Zha} for an elaborate explanation and application of least-squares methods in computer vision. Instead, we study parameter estimation by means of filtering. But before that, we mention a few merits and demerits of least-squares and other methods defined by algebraic distances. 

The use of algebraic distances in defining a cost function is computationally efficient and closed-form solutions are possible. The end result, however, is not satisfactory. This is due, in one part, to the fact that the objective function is mostly not invariant with respect to Euclidean transformations, for example, translation. This limits the coordinates systems to be used. In the other part, outliers may not contribute to the parameters the same way as inliers \cite{Zha}. Other more satisfactory parameter estimation methodologies are highly desired. We consider, in this paper, the use of Bayesian inference techniques.

Estimation of parameters by means of a filter can be achieved in a number of ways; one of them being the use of the filter evidence, or its near approximation, and selection criteria for parameters which give a reasonable estimate of the evidence. The second method involves updating the parameters and the state at the same time. This is known as dual estimation, which further subdivides into joint estimation and a dual filter. Joint estimation entails subjoining the vector of parameters to the state vector to form an extended state-space. The filter is then implemented and run forward in time with the hope of filter convergence to the optimal state and parameter values. A dual filter, on the other hand, involves implementing a filter for the state and that of parameters simultaneously. The filter provides a self-correcting mechanism which may lead to convergence of state and parameter estimates. 

Feedback particle filters \cite{Yan1,Yan2,Yan3, Yan4, Yan5, Ami, Amir} arose from the problem of particle degeneracy in particle filters, where particles lose weight substantially leading to poor filter performance. Feedback particle filters have been found to perform well and research is ongoing on the convergence and long-time performance of these filters \cite{Chenx, Path}. Other works related to this include time-delayed feedback control \cite{Xu, Zach}; the distinction from feedback particle filters being that the time-delay feedback control allows for direct control of parameters while the later uses a feedback gain which involves solving an associated partial differential equation \cite{Karl}.

The contribution of this paper is as follows: a new feedback particle filter is introduced based on the stochastic perturbation of the innovation process. this filter is shown to be exact. The feedback mechanism is then formulated from the perspective of the Schr{\"o}dinger problem. Sinkhorn algorithm is used to bridge the feedback mechanism hence solving, iteratively, the Schr{\"o}dinger problem. This yields the Sinkhorn particle filter. We investigate the performance of this filter in dual estimation and make comparisons with the performance of other extant filters. 

The balance of this paper is arranged as follows. We first introduce different filters in time-continuous framework. This is done in \cref{sec:filt}. The second part, that is \cref{sec:dual}, is devoted to dual estimation where extended state-space formulation and the dual filter mechanism are elaborated. Application of the foregoing ideas to estimating constant parameters forms the closing part of this paper. Experimental
results are in \cref{sec:experiments}, the discussion of results are in \cref{sec:discussions}, and the conclusions follow in
\cref{sec:conclusions}.

\section{Filters}
\label{sec:filt}
In this section, we introduce the filters used in this study.
\subsection{FPF with stochastically perturbed innovation} 
The feedback particle filter (FPF) \cite{Yan1,Yan2,Yan3, Yan4, Yan5, Ami, Amir} arose from the need to overcome the challenges that come with application of existing filters to nonlinear state-space models. Extended Kalman filters, although meant to exploit the optimal performance of the Kalman filters in nonlinear systems, are greatly hampered by the need to compute the Jacobians. This challenge is more prevalent in approximate second-order filters, what with the need to compute the Hessians besides. Moreover, the computational cost in the aforementioned filters is compounded by the need to integrate the equation for evolution of covariance. Whereas particle filters hold promise for overcoming the challenge faced in using approximate filters in nonlinear settings, they, too, have their own challenges. Particle filters are sequential sampling Monte Carlo methods with particle degeneracy mitigating resampling steps. Although they are easy to implement, they suffer from the so-called sample impoverishment owing to frequent resampling---which is why resampling is done when the effective sample size has deteriorated to a given threshold. Despite the improvisations tending to the improvement of particle filters, the particle filters retain a comparatively high error covariance.

The idea behind feedback particle filters is designing the posterior density of each particle so that it matches, optimally, the true posterior. As will be seen in the ensuing argument, if one begins by setting the prior density of each particle to be the same as the true prior density, FPF yields the same posterior density as the true posterior at all subsequent times. This matching of the posterior density of the particles with the true posterior is arrived at by minimizing a cost function. To the present time, the cost function used is the Kullback-Leibler divergence between the two posteriors. Minimization of the cost function, in the first variation, yields a Poisson equation, the solution of which results in an optimal control input. The FPF is characterized by a controlled innovation process and the gain. It turns out, and this will be apparent shortly, that the FPF gain is obtained by differentiating the optimal control input with respect to the particles.

Motivated by the two variants of Ensemble Kalman Bucy filter (EnKBF) based on the perturbation of the innovation process \cite{Law, Rei, Ber}, we obtain a stochastically-perturbed-innovation variant of the feedback particle filter by replacing the deterministic perturbation term to the innovation with a stochastic term. The McKean-Vlasov stochastic differential equations for the proposed feedback particle filter, for the state-space model \cref{eq:5.1.1a,eq:5.1.1b}, omitting the parameter term, consists of the equation of the mean-field process,
\begin{equation}
d\bar{x}_{t}  = f(\bar{x}_{t})dt +g(\bar{x}_{t}) d\bar{\beta}_{t} +  K(\bar{x}_{t})\circ (dy_{t} + R^{1/2}(t)d\bar{\eta}_{t}-h(\bar{x}_{t})dt), \label{eq:4.sd.1m}
\end{equation}
and the finite $M$ system of equations for the evolution of interacting hypotheses of the state, $X_{t}:= \{x^{i}_{t} \}^{M}_{i=1}$; that is,
\begin{equation}
dx^{i}_{t}  = f(x^{i}_{t})dt +g(x^{i}_{t}) d\beta^{i}_{t} +  K(x^{i}_{t})\circ (dy_{t} + R^{1/2}(t)d\eta^{i}_{t}-h(x^{i}_{t})dt), \label{eq:4.sd.1n}
\end{equation}
where $\bar{\beta}_{t}$ and $\{\beta_{t}\}_{i=1}^{M}$ are independent copies of $\beta_{t}$ and $\bar{\eta}_{t}$ and $\{\eta_{t}\}_{i=1}^{M}$ are independent copies of $\eta_{t}$; such that, for a given function $q(x)$, $$\pi_{t}[q(x)\mid Y_{t}] = \bar{\pi}_{t}[q(x)\mid Y_{t}] \approx \dfrac{1}{M} \sum_{i=1}^{M} q(x^{i}_{t}),$$
where $\pi_{t} = \textrm{Law}(x_{t})$ and $\pi_{t} = \textrm{Law}(\bar{x}_{t}).$

For notational convenience, we rewrite the \emph{controlled Stratonovich SDE} \cref{eq:4.sd.1n} as
\begin{equation}
dx_{t}  = f(x_{t})dt +g(x_{t}) d\beta_{t} +  K(x_{t})\circ (dy_{t} + R^{1/2}(t)d\eta_{t}-h(x_{t})dt), \label{eq:4.sd.1}
\end{equation}
where $K=\grad \phi(x)$ is the gain and $\phi(x)$ satisfies the Poisson equation, \begin{subequations}
\begin{align}
\div(\pi_{t}(x \mid Y_{t}) \grad \phi(x)) & = -(h(x_{t},t) - \hat{h})\pi_{t}(x \mid Y_{t})R^{-1}(t), \label{eq:4.1.3a}\\
\int \phi (x_{t})\pi_{t}(x \mid Y_{t})dx  & = 0, \label{eq:4.1.3b}
\end{align}
\end{subequations}
where $\pi_{t}(x \mid Y_{t})$ is the conditional density of the particle, $x^{i}_{t}$, given measurements $Y_{\tau}= \{y_{\tau} : 	t_{0} \leq \tau \leq t \} $. $\{\eta_{t}\}$ is an $r$-dimensional vector standard Brownian motion process.

Written in It\^{o} form, \cref{eq:4.sd.1} becomes
\begin{equation}
\begin{split}
dx_{t}  & = f(x_{t},t)dt +g(x_{t},t)d\beta_{t} + 2q(x_{t},t)dt \\
& \quad{} +  K(x_{t},t) (dy_{t} + R^{1/2}(t)d\eta_{t}-h(x_{t},t)dt), \label{eq:4.sd.2}
\end{split}
\end{equation}
where $q(x_{t},t)$ is given by 
\begin{equation}
 q_{j}(x_{t},t)= \dfrac{R}{2}\sum_{k=1}^{n}K_{k}(x,t)\dfrac{\partial K_{j}}{\partial x_{k}},\label{eq:4.1s.3}
\end{equation}

As in the FPF with a deterministically perturbed innovation, the equation of evolution of the filtering density $\pi (x_{t} \mid Y_{t})$ in the FPF with stochastically perturbed innovation is given by the Fokker-Planck equation to the SDE, \cref{eq:4.sd.2}; that is,
\begin{equation}
\begin{split}
d\pi_{t} & =  \mathcal{L}(\pi_{t} )dt - \div (\pi_{t}K)dy_{t} - \div (\pi_{t}U)dt  \\
& \quad{} + \sum_{k,l=1}^{n}\dfrac{\partial^{2}}{\partial x_{k}\partial x_{l}}(\pi_{t}(KRK^{\textrm{T}})_{kl})dt, \label{eq:4.sd.3}
\end{split}
\end{equation}
where $K=\grad \phi$ is the solution of 
\begin{equation}
\div (\pi_{t}K) = -\pi_{t} (h-\hat{h}_{t})R^{-1}(t) \label{eq:4.sd.4}
\end{equation}
and $U$ is then defined by
\begin{equation}
U = - Kh + 2q. \label{eq:4.sd.5}
\end{equation}

\subsubsection{Exactness of FPF with stochastically perturbed innovation} 

\begin{theorem}[Exactness of FPF with stochastic innovation]
Let initial filter posterior, $\pi_{t_{0}}$, be  equal to the true posterior $\pi^{*}_{t_{0}}$ at time $t_{0}$. Then the filter governed by \cref{eq:4.sd.1} is exact.  
\end{theorem}
\begin{proof}
We now show, that beginning with the filter posterior 
\begin{equation}
 \pi_{t_{0}}(x \mid Y_{t_{0}}) = \pi^{*}_{t_{0}}(x \mid Y_{t_{0}}), \label{eq:4.s.7}
 \end{equation} 
where $\pi^{*}_{t_{0}}(x \mid Y_{t_{0}})$ is the true posterior at initial time, $t_{0}$, the filter posterior matches the true posterior at all times, $t$; for then we demonstrate that the filter defined by \cref{eq:4.sd.1} is exact. It suffices to show that the equations of evolution of the true posterior and filter posterior are the same; that is, \cref{eq:4.sd.3} and 
\begin{equation}
d\pi^{*}_{t} =  \mathcal{L}(\pi^{*}_{t} )dt + \pi^{*}_{t} (h-\hat{h}_{t})^{\textrm{T}}R^{-1}(t)(dy_{t}-\hat{h}_{t}dt), \label{eq:4.1.4}
\end{equation} 
are the same. 

Multiplying $U$ in \cref{eq:4.sd.5} with $-\pi_{t}$ yields
\begin{equation}
\begin{split}
-\pi_{t}U &=  - \pi_{t}Kh + 2\pi_{t} q \\
& = - \pi_{t}K(h -\hat{h}_{t}) - \pi_{t}K\hat{h}_{t} + 2\pi_{t} q, \label{eq:4.s.8}
\end{split}
\end{equation}
where in the second equality we introduce $\pi_{t}K\hat{h}_{t} - \pi_{t}K\hat{h}_{t}$. From \cref{eq:4.sd.4} we have
\begin{equation}
-\pi_{t} (h-\hat{h}_{t}) = R(t) \div (\pi_{t}K_{k}). \label{eq:4.sd.9}
\end{equation}
Substituting \cref{eq:4.sd.9} in \cref{eq:4.s.8} gives
\begin{equation}
-\pi_{t}U  = KR(t) \div(\pi_{t}K)  -  \pi_{t}K\hat{h}_{t} + 2\pi_{t} q. \label{eq:4.s.9}
\end{equation}
Using 
\begin{equation}
\div\left( \pi \left[ KRK^{\textrm{T}}\right] \right) = \pi KR \div\left( K\right)  + KR \div\left( \pi K\right). \label{eq:4.3.3}
\end{equation} 
in \cref{eq:4.s.9} and noting the expression of $q(x_{t},t)$ in \cref{eq:4.1s.3}, we get
\begin{equation}
-\pi U =- \div\left( \pi \left[ KRK^{\textrm{T}}\right] \right)  + \pi K\hat{h}_{t}.  \label{eq:4.s.10}
\end{equation}
whereupon taking divergence on both sides of \cref{eq:4.s.10} yields
 \begin{equation}
 -\div( \pi U)  =- \sum_{l,k=1}^{n}\dfrac{\partial^{2}}{\partial x_{l} \partial x_{k}}\left( \pi \left[ KRK^{\textrm{T}}\right]_{lk} \right)   + \div(  \pi K\hat{h}_{t} ).  \label{eq:4.3s.4}
 \end{equation}
Finally, we substitute \cref{eq:4.sd.4,eq:4.3s.4} for $$-\div (\pi_{t}U) + \sum_{l,k=1}^{n}\dfrac{\partial^{2}}{\partial x_{l} \partial x_{k}}\left( \pi \left[ KRK^{\textrm{T}}\right]_{lk} \right)$$ and $\div (\pi_{t}K)$ in \cref{eq:4.sd.3} to obtain
\begin{equation}
d\pi_{t} = \mathcal{L}(\pi_{t} )dt + \pi_{t} (h-\hat{h}_{t})^{\textrm{T}}R^{-1}(t)(dy_{t}-\hat{h}_{t}dt),\label{eq:4.ss.3}
\end{equation}
which proves exactness.
\end{proof}

\begin{remark}[The gain]
The equation from which the gain is obtained is the same for both FPF with deterministically perturbed innovation and that with stochastically perturbed innovation.
\end{remark}

\subsection{Ensemble Transform Particle Filter} \label{sec:5.5}

The ensemble transform particle filter (ETPF) \cite{Rei1,Rei} is, as the FPF, characterised by a  control law, in the sense of optimal transportation, which moves the particles to the convenient position in the state space, thus improving upon the performance of the filter. The notable distinction between ETPF and particle filters with resampling is that the resampling step is replaced with a linear transformation. The transformation seeks to establish an optimal coupling between the prior ensemble and the posterior ensemble. 

Most precisely, the linear transport problem is as follows:
\begin{equation}
T^{*} = \underset{T\in \mathcal{R}^{M \times M}}{\textrm{arg min}} \sum_{i,j=1}^{M}t_{ij}\Vert x^{i}_{t_{n}}- x^{j}_{t_{n}} \Vert^{2}, \label{eq:6.4.12et}
\end{equation}
under the constraints
\begin{equation}
\sum_{i=1}^{M}t_{ij} = \dfrac{1}{M}, \; \sum_{j=1}^{M} t_{ij} = w_{t_{n}}^{i}, \; \textrm{  and  } t_{ij}>0,
\end{equation}
where $\{w^{i}_{t_{n}}\}_{i=1}^{M}$ are weights at time $t_{n}=n\delta t$ and $t_{ij}$ represents the element in row $i$ and column $j$ of the $M\times M$ matrix $T$. The weights are propagated as follows:
\begin{equation}
dw_{t}= \dfrac{1}{R}w_{t}^{i}(h(x^{i}_{t}) - \hat{h}_{t} )^{\textrm{T}}(dy_{t}-\hat{h}_{t}dt), \label{eq:5.4.14et}
\end{equation}
where \[\hat{h}_{t} = \dfrac{1}{M}\sum_{i=1}^{M}h(x^{i}_{t}).\]
To avoid negative weights, \cref{eq:5.4.14et} is approximated by
\begin{equation}
w^{i}_{t_{n+1}} \approx w^{i}_{t_{n}} \exp\left(-\dfrac{1}{2R}(h(x^{i}_{t_{n}}))^{T}h(x^{i}_{t_{n}})\delta t- (h(x^{i}_{t_{n}}))^{T}\delta y_{t_{n}} \right), \label{eq:6.4.15et}
\end{equation}
and initialized with \[w^{i}_{t_{0}}=\dfrac{1}{M}.\]

For the ETPF, the importance weights of the prior samples $\{x^{i}\}_{i=1}^{M}$ are given by $W_{1}= \{w_{1}^{i}\}_{i=1}^{M}$ and $W_{2}= \{w_{2}^{i}=1/M \}_{i=1}^{M}$. Moreover, $X_{1}=X_{2}=\{x^{i}\}_{i=1}^{M}$. The posterior samples are given by
\begin{equation}
\tilde{x}^{j}_{t_{n}} = \dfrac{1}{M}\sum_{i=1}^{M}x^{i}_{t_{n}}t_{ij}^{*}. \label{eq:6.4.19et}
\end{equation}

\subsection{Feedback formulation based on the Schr{\"o}dinger problem}

The Schr{\"o}dinger problem \cite{Schr} arises in the context of boundary value problems involving stochastic differential equations. Filtering problem can be cast as a Schr\"{o}dinger problem, by means of Bayes' rule, and in the space of probability measures \cite{Ang1,reich_2019}. In the ensuing argument, we invoke the algorithms used in optimal transportation \cite{Pey} to solve the Schr{\"o}dinger problem, which solution provides the estimates of the filtering distribution.

\begin{definition}{\textbf{Schr{\"o}dinger problem:}---} Find two functions $\hat{\phi}_{t_{n}}(x_{t_{n}})$ and $\phi_{t_{n+1}}(x_{t_{n+1}})$, satisfying the following equations
\begin{subequations}
\begin{align}
\pi_{t_{n}}(x_{t_{n}} \mid Y_{t_{n}}) & = \pi^{\phi}_{t_{n}}(x_{t_{n}} \mid Y_{t_{n}}) \hat{\phi}_{t_{n}}(x_{t_{n}}), \label{eq:4.6sh.1a} \\
\pi_{t_{n+1}}(x_{t_{n+1}} \mid Y_{t_{n+1}}) & = \pi^{\phi}_{t_{n+1}}(x_{t_{n+1}} \mid Y_{t_{n+1}}) \phi_{t_{n+1}}(x_{t_{n+1}}), \\
\pi^{\phi}_{t_{n+1}}(x_{t_{n+1}} \mid Y_{t_{n+1}}) & =  \int \pi_{t_{n+1}}(x_{t_{n+1}} \mid x_{t_{n}}) \pi^{\phi}_{t_{n}}(x_{t_{n}} \mid Y_{t_{n}})dx_{t_{n}}, \\
 \hat{\phi}_{t_{n}}(x_{t_{n}}) & = \int \pi_{t_{n+1}}(x_{t_{n+1}} \mid x_{t_{n}}) \phi_{t_{n+1}}(x_{t_{n+1}}) dx_{t_{n+1}}, \label{eq:4.6sh.1d} 
\end{align}
\end{subequations}
where $\pi_{t_{n}}(x_{t_{n}} \mid Y_{t_{n}})$ and $\pi_{t_{n+1}}(x_{t_{n+1}} \mid Y_{t_{n+1}})$ are the marginal filtering distributions at times $t_{n}$ and $t_{n+1}$, respectively.
\end{definition}
When solved, \crefrange{eq:4.6sh.1a}{eq:4.6sh.1d} yield transition density
\begin{equation}
\pi^{\phi}_{t_{n+1}}(x_{t_{n+1}} \mid x_{t_{n}}) := \dfrac{\phi_{t_{n+1}}(x_{t_{n+1}})}{\hat{\phi}_{t_{n}}(x_{t_{n}})}\pi_{t_{n+1}}(x_{t_{n+1}} \mid x_{t_{n}}), \label{eq:5.6.2}
\end{equation}
such that
\begin{equation}
\pi_{t_{n+1}}(x_{t_{n+1}} \mid Y_{t_{n+1}}) := \int \pi^{\phi}_{t_{n+1}}(x_{t_{n+1}} \mid x_{t_{n}})\pi_{t_{n}}(x_{t_{n}} \mid Y_{t_{n}})dx_{t_{n}}, \label{eq:5.6.3}
\end{equation}
where $\pi^{\phi}_{t_{n+1}}(x_{t_{n+1}}\mid x_{t_{n}})$ can be intuitively understood to define a coupling between the probability density function at time $t_{n}$, $\pi_{t_{n}}(x_{t_{n}} \mid Y_{t_{n}})$, and the filtering density function at time $t_{n+1}$, $\pi_{t_{n+1}}(x_{t_{n+1}} \mid Y_{t_{n+1}})$. 

How do these notions fit in to filtering? Recall that, in the discrete setting, Bayes' Theorem gives an expression of the filtering density, $\pi_{t_{n+1}}(x_{t_{n+1}} \mid Y_{t_{n+1}})$, at time $t_{n+1}$ as 
\begin{equation}
\pi_{t_{n+1}}(x_{t_{n+1}} \mid Y_{t_{n+1}}) = \dfrac{\pi_{t_{n+1}}(y_{t_{n+1}} \mid x_{t_{n+1}}) \pi_{t_{n+1}}(x_{t_{n+1}} \mid Y_{t_{n}})}{ \int \pi_{t_{n+1}}(y_{t_{n+1}} \mid x_{t_{n+1}}) \pi_{t_{n+1}}(x_{t_{n+1}} \mid Y_{t_{n}}) dx_{t_{n}}}, \label{eq:5.6.4}
\end{equation}
where 
\begin{equation}
\pi_{t_{n+1}}(x_{t_{n+1}} \mid Y_{t_{n}}) = \int \pi_{t_{n+1}}(x_{t_{n+1}} \mid x_{t_{n}})\pi_{t_{n}}(x_{t_{n}} \mid Y_{t_{n}})dx_{t_{n}} \label{eq:5.6.5}
\end{equation}
is the prediction density at time $t_{n+1}$ \cite{Sar,Jaz}. It then becomes apparent that \cref{eq:5.6.3,eq:5.6.4} are equivalent. Notice that the Schr{\"o}dinger problem reduces to \cref{eq:5.6.3} where, beginning with the filtering distribution $\pi_{t_{n}}(x_{t_{n}} \mid Y_{t_{n}})$, we directly obtain the filtering distribution $\pi_{t_{n+1}}(x_{t_{n+1}} \mid Y_{t_{n+1}})$ at time $t_{n+1}$ via the twisted transition density $\pi^{\phi}_{t_{n+1}}(x_{t_{n+1}} \mid x_{t_{n}})$.

We use empirical estimates of the densities to solve \crefrange{eq:4.6sh.1a}{eq:4.6sh.1d} governing the Schr\"{o}dinger problem; that is, given an ensemble $\{x^{i}_{t_{n}}\}_{i=1}^{M}$ at time $t_{n}$, we obtain an estimate of the prediction density, \cref{eq:5.6.5}, by 
\begin{equation}
\pi_{t_{n+1}}(x_{t_{n+1}} \mid Y_{t_{n}}) = \dfrac{1}{M} \sum_{i=1}^{M} \pi_{t_{n+1}}(x_{t_{n+1}} \mid x^{i}_{t_{n}}).
\end{equation}
The filtering density, \cref{eq:5.6.5}, is then empirically estimated as follows:
\begin{equation}
\pi_{t_{n+1}}(x_{t_{n+1}} \mid Y_{t_{n+1}}) = \dfrac{ \dfrac{1}{M} \sum_{i=1}^{M}\pi_{t_{n+1}}(y_{t_{n+1}} \mid x_{t_{n+1}}) \pi_{t_{n+1}}(x_{t_{n+1}} \mid x^{i}_{t_{n}})}{ \int \pi_{t_{n+1}}(y_{t_{n+1}} \mid x_{t_{n+1}}) \pi_{t_{n+1}}(x_{t_{n+1}} \mid Y_{t_{n}}) dx_{t_{n}}}.
\end{equation}
To solve the Schr{\"o}dinger problem, \crefrange{eq:4.6sh.1a}{eq:4.6sh.1d}, we make the following proposition
\begin{subequations}
\begin{align}
\pi^{\phi}_{t_{n}}(x_{t_{n}} \mid Y_{t_{n}}) & = \dfrac{1}{M} \sum_{i=1}^{M} \alpha^{i} \delta(x_{t_{n}}- x^{i}_{t_{n}}), \label{eq:4.6sh.4a} \\
\hat{\phi}_{t_{n}}(x_{t_{n}}) & = \dfrac{1}{M} \sum_{i=1}^{M} \dfrac{1}{\alpha^{i}}, \label{eq:4.6sh.4b}
\intertext{where}
\dfrac{1}{M} \sum_{i=1}^{M} \alpha^{i} & = 1. \label{eq:4.6sh.4c}
\end{align}
\end{subequations}
Substituting \cref{eq:4.6sh.4a,eq:4.6sh.4b} in to \crefrange{eq:4.6sh.1a}{eq:4.6sh.1d} yields,
\begin{subequations}
\begin{align}
\pi_{t_{n}}(x_{t_{n}} \mid Y_{t_{n}}) & = \dfrac{1}{M} \sum_{i=1}^{M} \delta(x_{t_{n}}- x^{i}_{t_{n}}), \label{eq:4.6sh.5a} \\
\pi_{t_{n+1}}^{\phi}(x_{t_{n+1}} \mid Y_{t_{n+1}})& = \dfrac{1}{M} \sum_{i=1}^{M} \alpha^{i} \pi_{t_{n+1}}(x_{t_{n+1}} \mid x_{t_{n}}^{i}), \label{eq:4.6sh.5b} \\
\phi_{t_{n+1}}(x_{t_{n+1}}) & = \dfrac{N}{ \int (N)dx_{t_{n}}\dfrac{1}{M} \sum_{i=1}^{M} \alpha^{i} \pi_{t_{n+1}}(x_{t_{n+1}} \mid x_{t_{n}}^{i})}, \label{eq:4.6sh.5c}
\end{align}
\end{subequations}
where $$N=\pi_{t_{n+1}}(y_{t_{n+1}} \mid x_{t_{n+1}}) \sum_{i=1}^{M}  \pi_{t_{n+1}}(x_{t_{n+1}} \mid x_{t_{n}}^{i}).$$ 
Therefore, one sets $\hat{\phi}_{t_{n}}(x_{t_{n}}) = \dfrac{1}{M} \sum_{i=1}^{M} \dfrac{1}{\alpha^{i}}$ and defines $\phi_{t_{n+1}}(x_{t_{n+1}})$ by \cref{eq:4.6sh.5c} in order to solve the Schr{\"o}dinger problem; for then, given a twisted prediction density
\begin{equation}
\begin{split}
\pi^{\phi}_{t_{n+1}}(x_{t_{n+1}} \mid Y_{t_{n}}) & = \dfrac{1}{M} \sum_{i=1}^{M} \pi^{\phi}_{t_{n+1}}(x_{t_{n+1}} \mid x^{i}_{t_{n}})\\
& =  \dfrac{1}{M} \sum_{i=1}^{M} \dfrac{\phi_{t_{n+1}}(x_{t_{n+1}})}{\hat{\phi}_{t_{n}}(x^{i}_{t_{n}})}\pi_{t_{n+1}}(x_{t_{n+1}} \mid x^{i}_{t_{n}}),
\end{split}
\end{equation}
by \cref{eq:5.6.2}, the relationship between the prediction density, $\pi^{\phi}_{t_{n+1}}(x_{t_{n+1}} \mid Y_{t_{n}}) $, and the twisted prediction density, $\pi^{\phi}_{t_{n+1}}(x_{t_{n+1}} \mid Y_{t_{n}}) $, is 
\begin{equation}
\dfrac{\pi_{t_{n+1}}(x_{t_{n+1}} \mid Y_{t_{n}}) }{\pi^{\phi}_{t_{n+1}}(x_{t_{n+1}} \mid Y_{t_{n}}) } = \dfrac{ \dfrac{1}{M} \sum_{i=1}^{M} \pi_{t_{n+1}}(x_{t_{n+1}} \mid x^{i}_{t_{n}})}{ \dfrac{1}{M} \sum_{i=1}^{M} \dfrac{\phi_{t_{n+1}}(x_{t_{n+1}})}{\hat{\phi}_{t_{n}}(x^{i}_{t_{n}})}\pi_{t_{n+1}}(x_{t_{n+1}} \mid x^{i}_{t_{n}})};
\end{equation}
from which, if we draw samples from the twisted distribution
\begin{equation}
x^{i}_{t_{n}} \sim \pi^{\phi}_{t_{n+1}}(x_{t_{n+1}} \mid x^{i}_{t_{n}}),
\end{equation}
then we can approximate the prediction density, $\pi^{\phi}_{t_{n+1}}(x_{t_{n+1}} \mid Y_{t_{n}}) $, by
\begin{subequations}
\begin{align}
\pi^{\phi}_{t_{n+1}}(x_{t_{n+1}} \mid Y_{t_{n}}) & \approx \dfrac{1}{M} w^{i} \delta (x_{t_{n}} - x_{t_{n}}), 
\intertext{where}
w^{i} & = \dfrac{\pi_{t_{n+1}}(x_{t_{n+1}} \mid Y_{t_{n}}) }{\pi^{\phi}_{t_{n+1}}(x_{t_{n+1}} \mid Y_{t_{n}}) }.
\end{align}
\end{subequations}

Now, suppose that we have $L=kM$, where $k\in \mathbb{N}$, samples $\{x_{t_{n}}^{j}\}_{j=1}^{L}$ drawn from $\pi_{t_{n}}(x_{t_{n}} \mid Y_{t_{n}})$. We obtain a bi-stochastic matrix $Q\in \mathbb{R}^{L\times M}$, which approximates the Markov process defined by $\pi(x_{t_{n+1}} \mid x_{t_{n}})$. The Schr{\"o}dinger problem, \crefrange{eq:4.6sh.1a}{eq:4.6sh.1d}, can then be recast as follows: Find two non-negative vectors $u\in \mathbb{R}^{L}$ and $v \in \mathbb{R}^{M}$ so that, 
\begin{equation}
P^{*} =  \textrm{diag}(u)Q \textrm{diag}(v)^{-1},
\end{equation}
given that $P^{*}$ belong to a polytope 
$$U := \left\lbrace  P \in \mathbb{R}^{L \times M} : P \geq 0, \sum_{j=1}^{L}p_{ji}= p_{1} , \;   \sum_{i=1}^{M}p_{ji} = p_{0}       \right\rbrace.$$
$P^{*}$ is also a solution to the optimization problem defined by minimizing the distance between all possible bi-stochastic matrices $P$ and $Q$; that is,
\begin{equation}
P^{*} = \underset{P \in U}{\textrm{arg min }} \textrm{KL}(P \Vert Q), \label{eq:4.6.6}
\end{equation}
where $\mathrm{KL}$ is the Kullback Leibler divergence between $P \in U$ and $Q$; that is,
\begin{equation}
\textrm{KL}(P \Vert Q) := \sum_{j,i =1}^{L,M} p_{ji} \log \dfrac{p_{ji}}{q_{ji}},
\end{equation}
where $p_{ji}$ and $q_{ji}$ are the elements of, respectively, matrices $P$ and $Q$ in row $j$ and column $i$.

Now the feedback particle filter in the Schr{\"o}dinger formulation is as follows: we first obtain an $M$-sized ensemble of states $\{ \bar{x}^{i}_{t_{n}}\}_{i=1}^{M}$ using a forecast distribution, $\pi_{t_{n}}( \bar{x}^{i}_{t_{n}} \mid Y_{t_{n}})$, with respect to which each particle evolves according to the weak form of the Fokker-Planck equation \cite{Ang1}; that is
\begin{equation}
\bar{x}^{i}_{t_{n}} = \bar{x}^{i}_{t_{n-1}} + f(\bar{x}^{i}_{t_{n-1}})\delta t. \label{eq:5.6.17}
\end{equation}
Secondly, $L=kM$ ensembles are obtained as follows: 
\begin{equation}
x^{j}_{t_{n}} = \bar{x}^{j}_{t_{n}} + g(t_{n})d\beta^{j}_{t_{n}}, \label{eq:5.6.18}
\end{equation}
where $\{\bar{x}^{j}_{t_{n}}\}_{j=1}^{L}$ are obtained by replicating each particle $\bar{x}^{i}_{t_{n}}$ $k$ times. The weights are obtained using
\begin{equation}
w^{j}_{t_{n}} = \exp \Bigl( -\dfrac{1}{2 \delta t} (- 2\delta y^{\textrm{T}}_{t_{n}} h(x^{j}_{t_{n}}) \delta t  + h^{\textrm{T}}(x^{j}_{t_{n}}) h(x^{j}_{t_{n}}) \delta t^{2} ) \Bigr) w^{j}_{t_{n-1}}. \label{eq:5.6.19}
\end{equation}
This gives a particle weight system, $\{x^{j}_{t_{n}}, \, \tilde{w}^{j}_{t_{n}}\}_{j=1}^{L}$, where $\tilde{w}^{j}_{t_{n}}$ signifies the normalised $w^{j}_{t_{n}}$---normalisation done according to $$\tilde{w}^{j}_{t_{n}} = \dfrac{w^{j}_{t_{n}}}{\sum_{j=1}^{L} w^{j}_{t_{n}}}.$$ 
By $Q \in \mathbb{R}^{L \times M}$, we denote and understand a matrix whose elements are
\begin{equation}
q_{ji} = \exp(-\dfrac{1}{2g^{2} \delta t} \Vert x^{j}_{t_{n}} -  \bar{x}^{i}_{t_{n}} \Vert^{2}  ). \label{eq:5.6.20}
\end{equation}
Then we solve the Schr\"{o}dinger problem defined by \cref{eq:4.6.6}, from whence we obtain $P^{*}$. $P^{*}$ can be obtained, iteratively, by means of the Sinkhorn scheme stipulated in the following algorithm, whose implementation details are stipulated in \cite{Cut}.
\begin{algorithm}[htp]
\caption{Sinkhorn iteration}
\label{alg:07}
\begin{algorithmic}[1]
\REQUIRE $p_{0}$ and $p_{1}$.  
\ENSURE $P$.
\STATE{Compute $Q$ and set $P^{0}=Q$}
\WHILE {$k>1$} 
\STATE{Compute: $u^{k+1} = \textrm{diag}(\sum_{m=1}^{M}P^{k})^{-1}p_{1}$. }
\STATE{Compute: $v^{k+1} = \textrm{diag}(\sum_{m=1}^{M}p_{0})^{-1} \sum_{l=1}^{L}(\textrm{diag}(u^{k+1})P^{k})^{\textrm{T}}$. }
\STATE{Compute: $P^{k+1} = \textrm{diag}(u^{k+1}) P^{k} \textrm{diag}(v^{k+1})^{-1}$. }
\ENDWHILE
\end{algorithmic} 
\end{algorithm}

Filtered particles are eventually obtained thus
\begin{equation}
 \tilde{x}^{i}_{t_{n}} = \sum_{j=1}^{L} x^{j}_{t_{n}} p_{ji}^{*} + g( \delta t)^{1/2}\xi^{i}_{t_{n}}, \qquad  \xi^{i}_{t_{n}} \sim \mathcal{N}(\bar{0}_{n \times 1}, \textrm{I}_{n\times n}), \label{eq:4.6.12}
\end{equation}
where $\bar{0}_{n \times 1}$ and $\textrm{I}_{n\times n}$ are, respectively, the null vector and the identity matrix of sizes as indicated in the subscripts. This yields the proposed Sinkhorn particle filter (SPF), whose summary is in \Cref{alg:05.6.2}.  

\begin{algorithm}[htp]
\caption{Sinkhorn particle filter}
\label{alg:05.6.2}
\begin{algorithmic}[1]
\REQUIRE $x^{i}_{t_{0}}$, $w^{i}_{t_{0}}= 1/M$ $\forall i \in \{1, \; 2, \; ..., \; M\}$, $k$, and $\delta y_{[t_{0},t_{T}]}$.  
\ENSURE $\hat{x}_{[t_{0},t_{T}]}$.
\FOR{$n=1$ \TO $N$, $\delta t>0$} 
\FOR{$i=1$ \TO $M$} 
\STATE{Obtain $\bar{x}^{i}_{t_{n}}$ using \cref{eq:5.6.17}  }
\FOR{$j=1$ \TO $L$} 
\STATE{Replicate $\bar{x}^{i}_{t_{n}}$ $k$ times and obtain $x^{j}_{t_{n}}$ using \cref{eq:5.6.17}  }
\STATE{Compute weights $w^{j}_{t_{n}}$ using \cref{eq:5.6.19}  }
\STATE{Compute $q_{ji}$ using \cref{eq:5.6.20}}
\ENDFOR
\STATE{Calculate $P^{*}$ by solving the the Schr{\"o}dinger problem via \Cref{alg:07} }
\STATE{Compute the filtered particles $\tilde{x}^{i}_{t_{n}}$ using \cref{eq:4.6.12} }
\ENDFOR
\STATE{Compute $\hat{x}_{t_{n}} = \dfrac{1}{M} \sum_{i=1}^{M} \tilde{x}^{i}_{t_{n}}$ }
\ENDFOR
\end{algorithmic} 
\end{algorithm}

Alternatively to obtaining the particles by means of \cref{eq:4.6.12}, we can resample $\{x^{j}_{t_{n}}\}_{j=1}^{L}$ such that
\begin{equation}
\mathbb{P} ( x^{i}_{t_{n}} = x^{j}_{t_{n}}  ) = p_{ji}^{*}, \qquad \forall i = 1, \, 2, \, ..., \, M. \label{eq:5.6.21}
\end{equation}
This results in the resampling Sinkhorn Particle filter (RSPF), the summary of which is in \Cref{alg:05.6.3}. 

\begin{algorithm}[htp]
\caption{Resampling Sinkhorn particle filter}
\label{alg:05.6.3}
\begin{algorithmic}[1]
\REQUIRE $x^{i}_{t_{0}}$, $w^{i}_{t_{0}}= 1/M$ $\forall i \in \{1, \; 2, \; ..., \; M\}$, $k$, and $\delta y_{[t_{0},t_{T}]}$.  
\ENSURE $\hat{x}_{[t_{0},t_{T}]}$.
\FOR{$n=1$ \TO $N$, $\delta t>0$} 
\FOR{$i=1$ \TO $M$} 
\STATE{Obtain $\bar{x}^{i}_{t_{n}}$ using \cref{eq:5.6.17}  }
\FOR{$j=1$ \TO $L$} 
\STATE{Replicate $\bar{x}^{i}_{t_{n}}$ $k$ times and obtain $x^{j}_{t_{n}}$ using \cref{eq:5.6.17}  }
\STATE{Compute weights $w^{j}_{t_{n}}$ using \cref{eq:5.6.19}  }
\STATE{Compute $q_{ji}$ using \cref{eq:5.6.20}}
\ENDFOR
\STATE{Calculate $P^{*}$ by solving the the Schr{\"o}dinger problem via \Cref{alg:07} }
\STATE{Resample the particles according to \cref{eq:5.6.21} to obtain the filtered particles $\tilde{x}^{i}_{t_{n}}$ }
\ENDFOR
\STATE{Compute $\hat{x}_{t_{n}} = \dfrac{1}{M} \sum_{i=1}^{M} \tilde{x}^{i}_{t_{n}}$ }
\ENDFOR
\end{algorithmic} 
\end{algorithm}

\section{Dual estimation}
\label{sec:dual}
Dual estimation comprehends simultaneous estimation of state and parameters by means of an appropriate filter. The self-correcting mechanism of the filter is taken advantage of to converge to both the true state and the true parameters. Depending on the initial parameter, the filter sooner or later converges to the true parameter value. Dual estimation can be achieved in two ways: joint estimation and by a dual filter \cite{Ang2,Lu,Lint,Mora,Ang}. In this section, and the rest of this paper, we shall assume that the parameters are static; that is, time-invariant.
\subsection{Joint estimation (augmented state-space)}
In joint estimation, the state vector is augmented with the vector of parameters to form an extended state-space and then the filter is run forward in time for an update of both the state and the parameters. The parameters are induced with artificial dynamics, or are made to assume a random walk; that is, respectively,
\begin{equation}
 dz_{t} = \zeta_{t} ; \quad t_{0} \leq t, \label{eq:6.4.1}
\end{equation}
where 
\begin{subequations}
\begin{align}
dz_{t} & =\begin{pmatrix}
dx_{t} \\
d\theta_{t}
\end{pmatrix} \textrm{  and  } \zeta_{t}=\begin{pmatrix}
f(x_{t},\theta)dt  + g(x_{t})d\beta_{t} \\
0
\end{pmatrix},
\intertext{or where}
\zeta_{t} & =\begin{pmatrix}
f(x_{t},\theta)dt  + g(x_{t})d\beta_{t} \\
\sigma d\chi_{t}
\end{pmatrix},
\end{align}
\end{subequations}
in which $\{\chi_{t},\, t>t_{0}\}$ is a $d-$dimensional standard Brownian motion vector process and $\sigma$ is a small constant. A filter is then implemented with the augmented state $z_{t}$ in the place of $x_{t}$. The demerit of this method is that the extended state-space has an increased degree of freedom owing to many unknowns, of both the state and the parameters, which renders the filter unstable and intractable, especially in nonlinear models \cite{Mor}. Joint estimation is also known as expectation maximisation \cite{Sar}.

\subsection{Test example}
We now consider a scaler example in order to test the performance of different filters under joint estimation of both the state and the parameters. 

\begin{exam}[Scalar SDE]{ex5}
\vspace{-0.4cm}
Consider the following linear Gaussian It\^{o} state space model.
\begin{subequations}
\begin{align}
dx_{t} &= (ax_{t} + b)dt + Q^{1/2}dv_{t}; \qquad  t_{0} \leq t, \\
dy_{t} &= cx_{t}dt + R^{1/2}dw_{t}; \qquad  t_{0} \leq t,
\end{align}
\end{subequations}
where $\{v_{t}\}$ and $\{w_{t}\}$ are standard Brownian motion processes with, respectively, \[ \mathbb{E}\{ dv_{t}dv_{t}^{T}\}=dt \textrm{  and  } \mathbb{E}\{ dw_{t}dw_{t}^{T}\}=dt.\] Let the state, $x_{t}$, at time $t_{0}$ be $x_{t_{0}}\sim \mathcal{N}(0,0.001)$. Let, moreover, $x_{t_{0}}$, $\{v_{t}, \, t \geq t_{0}\}$ and $\{v_{t}, \, t \geq t_{0}\}$ be uncorrelated. We take $ a=-0.2$, $b=0.2$, $c=1.01$, $Q=0.001$, $R=0.001$ and proceed to simultaneously estimate the state, $x_{t}$, and the parameters, $\textcolor{red}{a}$ and $\textcolor{red}{b}$ using different filters. 
\end{exam}

\section{Experimental results}
\label{sec:experiments}

The following panels show the results obtained using EnKBF, BPF (bootstrap particle filter) \cite{Ang1,Aru}, FPF (feedback particle filter with kernel-based gain approximation), ETPF---ensemble transform particle filter \cite{Rei} introduced in \Cref{sec:5.5} and RSPF. 

\begin{figure}[htp]
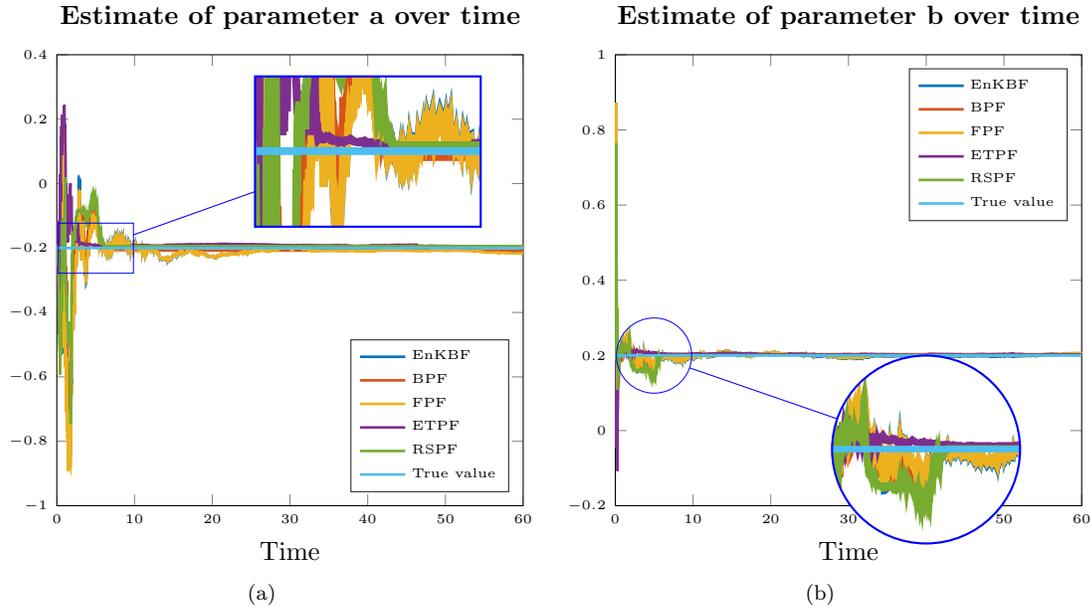

\centering
\subfigure[Estimates of parameter $a$][]{
\input{Images/Pest5.tex}
\label{fig:subfig 5.2c}}
\subfigure[Estimates of parameter $b$][ ]{
\input{Images/Pest6.tex}
\label{fig:subfig 5.2d}}
\caption[Plots showing:- (a) estimates of parameter $a$ and (b) estimates of parameter $b$ over time using EnKBF, BPF, FPF, ETPF and RSPF]{Plots showing:- (a) estimates of parameter $a$ and (b) estimates of parameter $b$ over time using EnKBF, BPF, FPF, ETPF and RSPF. The true parameter values are, respectively, $a=-0.2$ and $b=0.2$. The plots indicate that all the filters converge to the true parameter values. The time step used is $\delta t=0.02$.}
\label{fig:globfig5.2b}
\end{figure}

We now plot the box-plots, the better to see the distribution of parameter estimates in the results shown in \Cref{fig:globfig5.2b}, beyond time $30$. 

\begin{figure}[htp]
\centering
\subfigure[Estimates of parameters $a$][]{
\includegraphics[width=0.45\textwidth]{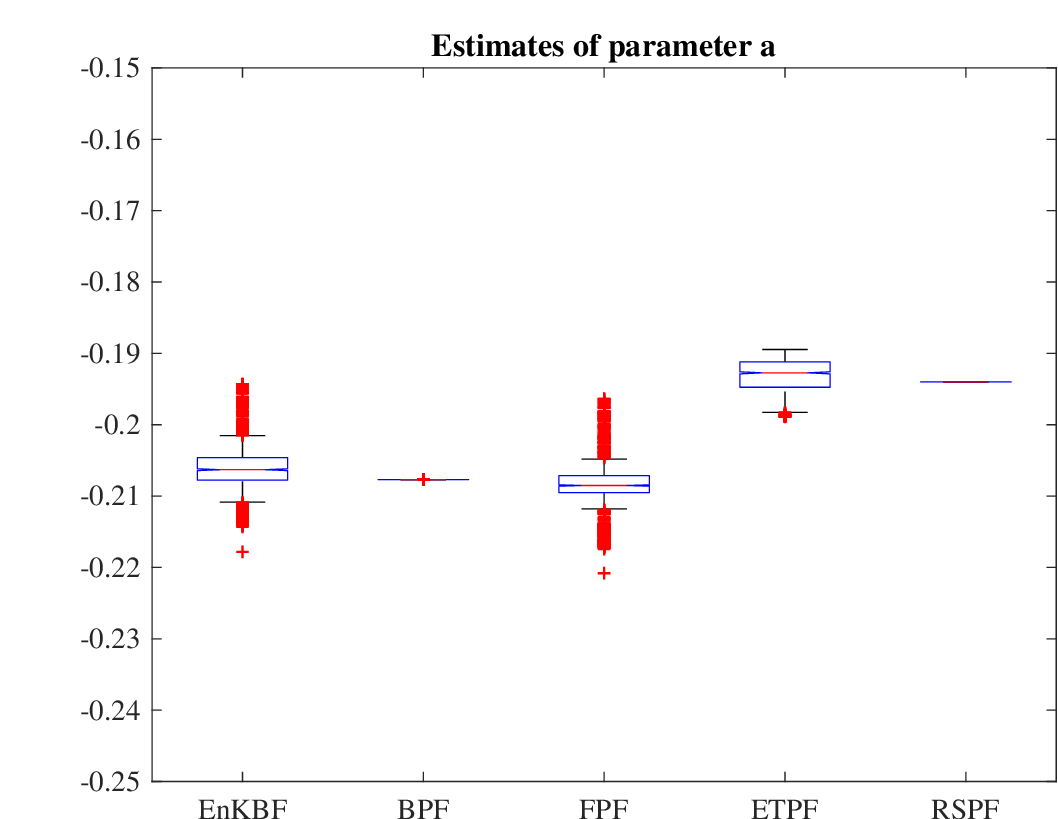}
\label{fig:subfig 5.2a}}
\subfigure[Estimates of parameters $b$][ ]{
\includegraphics[width=0.45\textwidth]{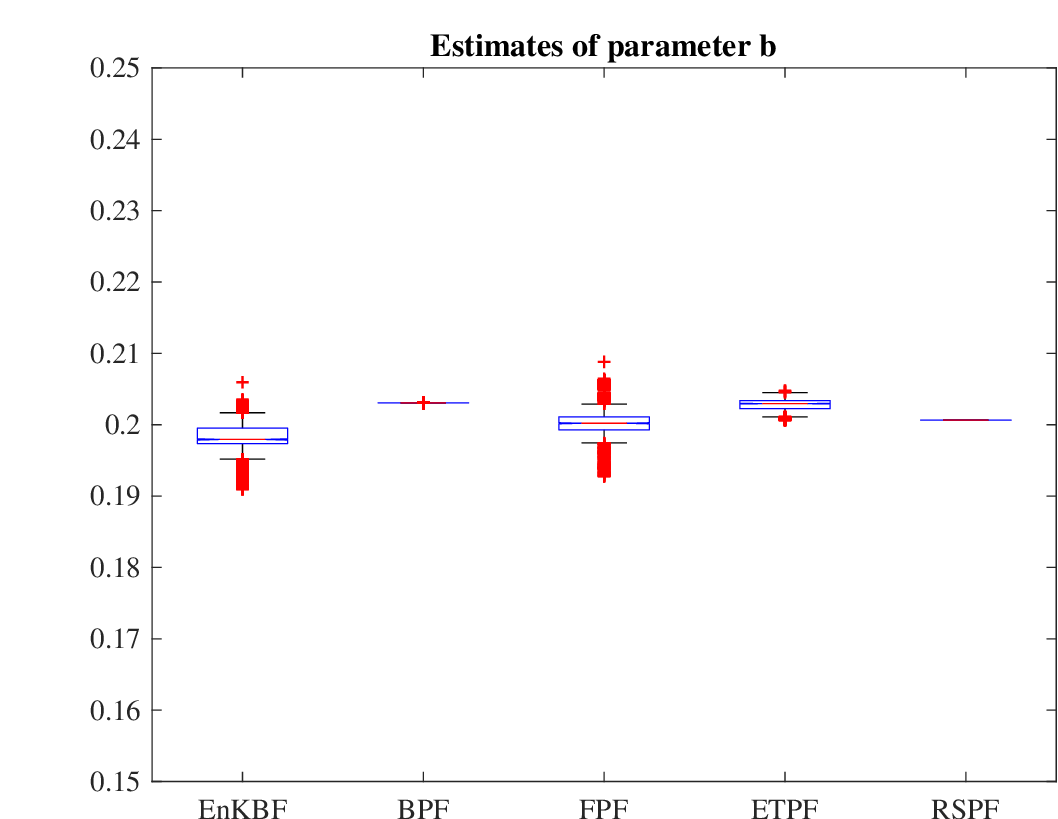}
\label{fig:subfig 5.2b}}
\caption[Box-plots showing the distribution of estimates of parameters $a$ and $b$ obtained using EnKBF, BPF, FPF, ETPF and RSPF]{Box-plots (a) and (b) showing, respectively, the distribution of estimates of parameters $a$ and $b$ obtained using EnKBF, BPF, FPF, ETPF and RSPF. The gain in the FPF is computed using the kernel based gain approximation method \cite{Ang1}. In both cases, the EnKBF and FPF register more dispersive results with more outliers than the BPF, ETPF and RSPF.}
\label{fig:globfig5.2}
\end{figure}

\begin{table}[htp]
\begin{center}
\begin{tabular}{llllll} 
\hline
&\textbf{EnKBF} & \textbf{BPF} & \textbf{FPF} & \textbf{ETPF} & \textbf{RSPF} \\ 
\hline
$M$ & $1000$& $1000$& $1000$& $100$ & $1000$\\ 
Time in seconds & $21.69$& $11.78$& $233.66$& $334.78$ & $732.12$ \\ 
\hline
\end{tabular} 
\caption[A table showing run-time and number of particles for results shown in \Cref{fig:globfig5.2} for different filters]{A table showing run-time and number of particles for results shown in \cref{fig:globfig5.2} for different filters. This is an output of 2.3 GHz Intel Core i5 processor. \label{tab:6.1}} 
\end{center}
 \end{table}
 
 \begin{center}
 \begin{figure}
 \begin{tikzpicture}
  \begin{axis}[title  = Root Mean Square Error (RMSE) of Different Filters,
    xbar,
    y axis line style = { opacity = 0 },
    ylabel={Filters},
    axis x line       = none,
    tickwidth         = 0pt,
    ytick             = data,
    enlarge y limits  = 0.2,
    enlarge x limits  = 0.02,
    symbolic y coords = {EnKBF, BPF, FPF, ETPF, RSPF},
    nodes near coords,
  ]
  \addplot coordinates { (0.010117972210033,EnKBF)         (0.006825200295576,BPF)
                         (0.013445477969967,FPF)  (0.0021293026630727,ETPF)  (0.002031136630727,RSPF)};
  \addplot coordinates { (0.003829670746118,EnKBF)         (0.002914768292107,BPF)
                         (0.005559297377356,FPF)   (0.0017435639323763,ETPF) (0.001537779323763,RSPF)  };
  \legend{$a$, $b$}
  \end{axis}
\end{tikzpicture}
\caption{Plot showing the RMSE in estimation of parameters a and b using EnKBF, BPF, FPF, ETPF and RSPF. ETPF and RSPF register a small RMSE which shows that they perform better than the rest of the filters. RSPF, however, although having the lowest RMSE, is run with 1000 ensembles while ETPF is run with 100 ensembles because of time constraints. }
\label{fig:globfig3.11}
\end{figure}
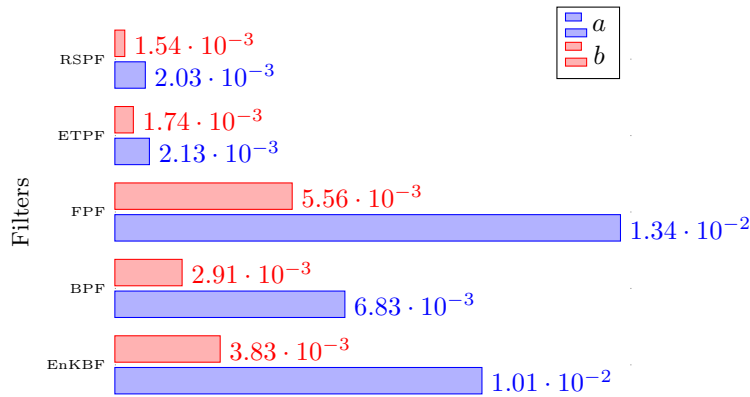
\end{center}

\section{Discussion of experimental results}
\label{sec:discussions}
The EnKBF, BPF, FPF, ETPF and RSPF yield converging results to the true parameter values but with some margin of error. From \Cref{fig:subfig 5.2c,fig:subfig 5.2d}, it is evident that the ETPF yields a faster converging results compared to the other filters. Furthermore, the FPF takes more time before a steady value is obtained as compared to other filters.

From \Cref{tab:6.1}, the BPF is faster compared to other methods. FPF is slower whilst ETPF is the slowest, seeing that the number of particles used is $100$. The gain in the FPF is computed using the kernel based gain approximation method. The reason for this is that FPF involves a lengthy procedure for gain computation, at each iteration. On the other hand, The ETPF involves frequent solution of optimal transport problem. In this example, earth movers distance (EMD) algorithms were used in the ETPF. As a remedy, faster optimal transport algorithms, for example the Sinkhorn iteration, can be viable \cite{Cut}. Although the RSPF takes much time, the result, as seen in \Cref{fig:globfig5.2}, is more accurate compared to that obtained from other methods. This is further shown in \Cref{fig:globfig3.11}, where the lowest RMSE is posted by RSPF indicating that it is more accurate compared to the rest of the filters.  

\section{Conclusions}
\label{sec:conclusions}
This paper has accomplished a number of purposes. In the first place, a variant of feedback particle filter has been introduced with stochastically perturbed innovation process and it has been shown that the proposed filter is exact; that is to say, the filter posterior maps the true posterior when the filter prior is set equal to the true prior. Sinkhorn particle filters have been introduced beginning with optimal transport notions as postulated by Schr{\"o}dinger. Finally, Comparison of the performance of EnKBF, BPF, FPF, ETPF and RSPF filters in joint estimation of both the state and parameters has been carried out by means of a scalar example. From the results, ETPF seems to converge to the true parameter faster that the other filters. In as much as RSPF takes more time, the boxplots indicate that RSPF is more accurate compared to EnKBF, BPF, FPF,and ETPF. Moreover, the plot of RMSE indicates that RSPF is more accurate compared with the other filters.

\section*{Acknowledgments}
I would like to acknowledge the assistance of Professor Dr. Sebastian Reich, of Potsdam University, Germany, in supervising my doctoral research \cite{Ang1}, of which this paper is part.

\bibliographystyle{siamplain}
\bibliography{mybib}

\begin{thebibliography}{10}

\bibitem{Ami}
{\sc T.~Amirhossein, J.~de~Wiljes, G.~Prashant, and S.~Reich}, {\em {Kalman
  Filter and its Modern Extensions for the Continuous-time Nonlinear Filtering
  Problem}}, Psychological Review,  (2016).

\bibitem{Ang1}
{\sc D.~Angwenyi}, {\em Time-continuous state and parameter estimation with
  application to hyperbolic SPDEs}, doctoralthesis, Universit{\"a}t Potsdam,
  2019, \url{https://doi.org/10.25932/publishup-43654}.

\bibitem{Ang2}
{\sc D.~Angwenyi}, {\em {Estimation of spatially varying parameters with
  application to hyperbolic SPDEs}}, arXiv:2107.07246v1,  (2021).

\bibitem{Ang}
{\sc D.~Angwenyi, J.~de~Wiljes, and S.~Reich}, {\em {Interacting Particle
  Filters for Simultaneous State and Parameter Estimation}}, arXiv:
  1709.09199v1,  (2017).

\bibitem{Aru}
{\sc S.~Arulampalam, S.~Maskell, N.~Gordon, and T.~Clapp}, {\em {A Tutorial on
  Particle Filters for Online Nonlinear/Non-Gaussian Bayesian Tracking}}, IEEE
  Transactions on Signal Processing, 50 (2002), pp.~174 -- 188.

\bibitem{Ber}
{\sc K.~Bergemann and S.~Reich}, {\em {An Ensemble Kalman-Bucy Filter for
  Continuous Data Assimilation}}, Meteorolog. Zeitschrift, 21 (2012), pp.~213
  -- 219.

\bibitem{Karl}
{\sc K.~Berntorp}, {\em {Feedback particle filter: Application and
  evaluation.}}, in 2015 18th International Conference on Information Fusion,
  2015, pp.~1633--1640.

\bibitem{Chenx}
{\sc X.~Chen, X.~Luo, J.~Shi, and S.~Yau}, {\em {General convergence result for
  continuous-discrete feedback particle filter.}}, International Journal of
  Control,  (2021).

\bibitem{Cut}
{\sc M.~Cuturi}, {\em {Sinkhorn Distances: Lightspeed Computation of Optimal
  Transportation Distances}}, arXiv:1306.0895v1 [stat.ML],  (2013).

\bibitem{Ein}
{\sc G.~A. Einicke}, {\em {Smoothing, Filtering and Prediction: Estimating the
  Past, Present and Future}}, InTech, Janeza Trdine 9, 51000 Rijeka, Croatia,
  2012.

\bibitem{Jaz}
{\sc A.~Jazwinski}, {\em {Stochastic Processes and Filtering Theory}}, Academic
  Press, New York, 1970.

\bibitem{Law}
{\sc K.~Law, A.~Stuart, and K.~Zygalakis}, {\em {Data Assimilation: A
  Mathematical Introduction}}, Springer, 2015.

\bibitem{Lew}
{\sc J.~Lewis, S.~Lakshmivarahan, and S.~Dhall}, {\em {Dynamic Data
  Assimilation: A Least Squares Approach}}, Cambridge University Press, 2006.

\bibitem{Lint}
{\sc J.~W. C.~V. Lint, S.~P. Hoogendoorn, and A.~Hagyi}, {\em {Dual EKF State
  and Parameter Estimation in Multi-Class First-order Traffic Models}}, in
  Proceedings of the17th World Congress, The International Federation of
  Automatic Control, Seoul, Korea, July 6-11 2008.

\bibitem{Lu}
{\sc H.~L\"{u}, Z.~Yu, Y.~Zhu, S.~Drake, Z.~Hao, and E.~A. Sudicky}, {\em {Dual
  State-Parameter Estimation of Root Zone Soil Moisture by Optimal Parameter
  Estimation and Extended Kalman Filter Data Assimilation}}, Advances in Water
  Resources, 34 (2011), pp.~395--406.

\bibitem{Mor}
{\sc H.~Moradkhani, S.~Sorooshian, H.~Gupta, and P.~Houser}, {\em {Dual
  State-Parameter Estimation of Hydrological Models Using Ensemble Kalman
  Filter}}, Advances in Water Resources, 28 (2005), pp.~135--147.
\newblock Elsevier.

\bibitem{Mora}
{\sc H.~Moradkhani, S.~Sorooshian, H.~V. Gupta, and P.~R. Houser}, {\em {Dual
  State-Parameter Estimation of Hydrological Models Using Ensemble Kalman
  Filter}}, Advances in Water Resources, 28 (2005), pp.~135--147.

\bibitem{Path}
{\sc S.~Pathiraja and W.~Stannat}, {\em {Analysis of the feedback particle
  filter with diffusion map based approximation of the gain.}}, Foundations of
  Data Science, 3 (2021).

\bibitem{Pey}
{\sc G.~Peyr\'{e} and M.~Cuturi}, {\em {Computational Optimal Transport}},
  arXiv: 1803.00567v1 [stat.ML], 2018.

\bibitem{Rei1}
{\sc S.~Reich}, {\em {A Nonparametric Ensemble Transform Method for Bayesian
  Inference}}, SIAM Journal of Scientific Computing, 35 (2013),
  pp.~A2013--A2024.

\bibitem{reich_2019}
{\sc S.~Reich}, {\em {Data assimilation: The Schrödinger perspective}}, Acta
  Numerica, 28 (2019), p.~635–711,
  \url{https://doi.org/10.1017/S0962492919000011}.

\bibitem{Rei}
{\sc S.~Reich and C.~Cotter}, {\em {Probabilistic Forecasting and Bayesian Data
  Assimilation}}, Cambridge University Press, 2015.

\bibitem{Sar}
{\sc S.~S\"{a}rkk\"{a}}, {\em {Bayesian Filtering and Smoothing}}, Cambridge
  University Press, Cambridge, 2013.

\bibitem{Schr}
{\sc E.~Schr\"{o}dinger}, {\em {\"{U}ber die Umkehrung der Naturgesetze}},
  Sitzungberichte der Preu{\upshape{\ss}}ischen Akademie der Wissenschaften,
  Physikalisch-mathematische Klasse,  (1931), pp.~144--153.

\bibitem{Amir}
{\sc A.~Taghvaei and P.~G. Mehta}, {\em Error analysis of the stochastic linear
  feedback particle filter}, arXiv: 1809.07892v1,  (2018).

\bibitem{Xu}
{\sc E.~Xu, Y.~Zhang, and Y.~Chen}, {\em {Time-delayed local feedback control
  for a chaotic finance system.}}, Journal of Inequalities and Applications,
  (2021).

\bibitem{Yan5}
{\sc T.~Yang}, {\em {Feedback Particle Filter and its Applications}}, PhD
  thesis, Department of Mechanical Science and Engineering, University of
  Illinois at Urbana-Champaign, 2014.

\bibitem{Yan3}
{\sc T.~Yang, R.~Laugesen, P.~Mehta, and S.~Meyn}, {\em {Multivariable feedback
  particle filter}}, in Proceedings of the 51st IEEE Conference on Decision
  Control, 2012, pp.~4063--4070.
\newblock Maui, HI.

\bibitem{Yan1}
{\sc T.~Yang, P.~Mehta, and S.~Meyn}, {\em {A Mean-field Control-oriented
  Approach to Particle Filtering}}, in Proceedings of American Control
  Conference, 2011, pp.~2037--2043.

\bibitem{Yan2}
{\sc T.~Yang, P.~Mehta, and S.~Meyn}, {\em {Feedback Particle Filter with
  Mean-field Coupling}}, in Proceedings of IEEE 50th Conference on Decision and
  Control, 2011, pp.~7909--7916.
\newblock Orlando, FL.

\bibitem{Yan4}
{\sc T.~Yang, P.~Mehta, and S.~Meyn}, {\em {Feedback Particle Filter}}, IEEE
  Transactions in Automatic Control, 58 (2013), pp.~2465--2480.

\bibitem{Zach}
{\sc A.~Zacharova, N.~Semenova, V.~Anishchenko, and E.~Schoell}, {\em
  {Time-delayed feedback control of coherence resonance chimeras.}}, An
  Interdisciplinary Journal of Nonlinear Science., 27 (2017).

\bibitem{Zha}
{\sc Z.~Zhang}, {\em {Parameter Estimation Techniques: A Tutorial with
  Applications to Conic Fitting}}, Image and Vision Computing,  (1997).
\newblock Elsevier.

\end{thebibliography}

\end{document}